\documentclass[9pt,draft]{amsart}

\usepackage[latin1]{inputenc}
\usepackage{latexsym}
\usepackage{a4wide}
\usepackage{amscd}
\usepackage{graphics}
\usepackage{amsmath}
\usepackage{amssymb}
\usepackage{mathrsfs}
\input xy
\xyoption{all}
\usepackage{bm}

\newcommand{\C}{\mathbb C}
\newcommand{\N}{\mathbb N}
\newcommand{\R}{\mathbb R}

\def\virgA{``}

\newcommand{\sgn}{\mathrm{sign}}

\newcommand{\ispec}{\mu_{\scriptstyle{\mathrm{Mor}}}}

\newcommand{\noo}[1]{\overset {\mbox{%
\lower1pt\hbox{${\scriptscriptstyle o}$}}}n^{\mbox{%
\lower2pt\hbox{$\scriptscriptstyle #1$}}}}

\newcommand{\la}{\lambda}

\newcommand{\spfl}{\mathrm{sf}}

\newcommand{\icon}{\mu_{\scriptstyle{\mathrm{con}}}}

\newcommand{\Id}{\mathrm{Id}}

\newcommand{\Imm}{\mathrm{Im}}

\newcommand{\diag}{\mathrm{diag\,}}

\newcommand{\finedim}{\hfill $\Box$}            
\numberwithin{equation}{section}
 \usepackage{geometry}

\newtheorem{mainthm}{\sc Theorem}           
\newtheorem{thm}{\sc Theorem}[section]      
\newtheorem{lem}[thm]{\sc Lemma}            
\newtheorem{defin}[thm]{\sc Definition}      

\title{Brief communication. An indefinite Sturm theory.}
\author[A. Portaluri]{Alessandro Portaluri}
\address{Dipartimento di Matematica \virgA Ennio De Giorgi ",\hfill\break\indent
Collegio Fiorini,  Università del Salento, \hfill\break\indent Via
per Arnesano, Lecce, Le\hfill\break\indent Italy}
\email{alessandro.portaluri@unile.it}
\date{December 10, 2008. Revised version.}
\thanks{%
The author was partially supported by MIUR project {\em Variational
Methods and Nonlinear Differential Equations\/}.}

\begin{document}

\maketitle

\section*{Introduction} Sturm theory for second order differential
equations was generalized to systems and higher order equations with
positive leading coefficient by several authors. (See for instance
\cite{Arn86}, \cite{Edw64} and \cite{Ovs90} among the others). Here
we propose a Sturm oscillation theorem for indefinite systems with
Dirichlet boundary conditions of the form
\begin{equation}\label{eq:opdiff}
l(x, D)\,v:= p_{2m}\dfrac{d^{2m}v}{dx^{2m}}+
p_{2m-2}(x)\dfrac{d^{2m-2}v}{dx^{2m-2}}+ \dots
+p_1(x)\dfrac{dv}{dx}+ p_0(x) v =0,
\end{equation}
where  $p_i$ is a smooth path of matrices on the complex
$n$-dimensional vector space $\C^n$, $p_{2m}$ is the symmetry of the
form $\diag(I_{n-\nu}, -I_\nu)$ for some $\nu \geq0$. This
generalization was obtained along the lines of \cite{MusPejPor03}
and \cite{MusPejPor04} in which the second order equations were
considered. Full proofs and the relation with Maslov index will
appear elsewhere.

\section{Variational set up and main results}

We will use the variational approach to \eqref{eq:opdiff} as
described in \cite{Edw64} and we will stick to the notations of that
paper. Given the complex $n$-dimensional Hermitian space $(\C^n,
\langle \cdot, \cdot \rangle)$, for any $m \in \N$ let $\mathscr
H^m:=H^m(J,\C^n)$ be the
Sobolev space of all $H^m$-maps from $J:=[0,1]$ into $\C^n$.\\
\noindent A {\em derivative dependent Hermitian form\/} is the form
$ \Omega(x)[u] = \sum_{i,j=0}^m \langle D^i u(x), \omega_{i,j}(x)
D^j u(x)\rangle $, where, each $\omega_{i,j}$ is a smooth path of
$x$-dependent Hermitian matrices with constant leading coefficient
$\omega_{m,m}:=p_{2m}$ and such that $\omega_{i, 2m-1-i}=0$ for each
$i=0, \dots, m$. Each derivative dependent Hermitian form, defines a
Hermitian form $q \colon \mathscr H^m \to \R$ by setting $ q(u):=
\int_J \Omega(x)[u] dx.$ If $v \in \mathscr H^m$ and $u \in \mathscr
H^{2m}$ then, using integration by parts, the corresponding
sesquilinear form $q(u,v)$ can be written as
\begin{equation}\label{eq:lachiama2.0}
q(u,v)= \int_J\langle v(x), l(x,D)u(x) \rangle dx + \phi(u,v),
\end{equation}
where $l(x,D)$ is a  differential operator of the form of
\eqref{eq:opdiff} and $\phi(u,v)$ is a bilinear form depending only
on the $(m-1)$-jet, $j^{m-1} v(x):=(v(x), \dots, v^{(m-1)}(x))$ and
on the $(2m-1)$-jet $j^{2m}u(x)$ at the boundary $x=0,1$.\\
\noindent Let $\mathscr H_0^m:=\mathscr H_0^m(J):=\{u \in \mathscr
H^m \colon j^{m-1} u(0)=0=j^{m-1} u(1)\}$ and let $q_\Omega$ be the
restriction of the Hermitian form $q$ to $\mathscr H_0^m$. For each
$\lambda \in J$, let us consider the space $\mathscr H_0^m([0,
\lambda])$ with the form $\int_{[0,\lambda]}\Omega(x)dx$. Via the
substitution $x \mapsto \lambda x$, we transfer this form to
$\mathscr H_0^m(J)$, so, we come to the forms $\Omega_\lambda$ and
$q_\lambda$  defined respectively by
\begin{equation}\label{eq:5}
\Omega_\lambda(x)[u] := \sum_{i,j=0}^{m} \langle D^i u(x),
\lambda^{2m-(i+j)}\omega_{i,j}(\lambda x) D^j u(x)\rangle \quad
\textrm{and we let}\quad q_\lambda(u):= \int_J
\Omega_\lambda(x)[u]dx.
\end{equation}
Then $\lambda \mapsto q_\lambda$ is a smooth path of Hermitian forms
acting on $\mathscr H_0^m$ with $q_1=q_{\Omega}$ and with $q_0(u)=
\int_J \langle p_{2m} D^m u, D^m u \rangle dx$. Now we introduce the
following definition.
\begin{defin}
A {\em conjugate instant\/} for $q_{\Omega}$ is any point $\lambda
\in(0,1]$ such that $\ker q_\lambda \not=\{0\}$.
\end{defin}

\noindent Let $C_\lambda$ be the path of bounded selfadjoint
Fredholm operators associated to $q_\lambda$  via the Riesz
representation theorem.
\begin{lem}
The Hermitian form $q_0$ is non-degenerate. Moreover each
$q_\lambda$ is a Fredholm Hermitian form. (i.e. $C_\lambda$ is a
Fredholm operator). In particular $\dim \ker q_\lambda <+\infty$ and
$q_\lambda$ is non degenerate if and only if $\ker q_\lambda =
\{0\}$.
\end{lem}
\proof That the operator $C_0$ is an isomorphism can be proven
exactly as in \cite[Proposition 3.1]{MusPejPor03}. On the other hand
each $q_\lambda$ is a weakly continuous perturbation of $q_0$ since
it differs from $q_0$ only by derivatives of $u$ of order less than
$m$. Therefore $C_\lambda- C_0$ is compact for all $\lambda \in J$
and hence $C_\lambda$ is Fredholm of index $0$. The last assertion
follows from this.\finedim

When the form $q_\Omega$ is non degenerate, i.e. when $1$ is not a
conjugate instant, we introduce the following definition.
\begin{defin}
The {\em (regularized) Morse index\/} of $q_\Omega$ is defined by
\begin{equation}\label{def:ispettrale}
\ispec(q_\Omega):= -\spfl (q_\lambda, J),
\end{equation}
where $\spfl$ denotes the {\em spectral flow of the path
$q_\lambda$\/} i.e., the number of positive eigenvalues of
$C_\lambda$ at $\lambda=0$ which become negative at $\lambda=1$
minus the number of negative eigenvalues of $C_\lambda$ which become
positive. (See, for instance \cite{MusPejPor03}, for a more detailed
exposition).
\end{defin}
Now we define the conjugate index. By \eqref{eq:lachiama2.0} $q(u,v)
= \int_J \langle v(x), l(x,D) u(x) \rangle dx$ for $u \in \mathscr
H_0^{2m}$ and $v \in \mathscr H_0^m$ and therefore elements of the
kernel of $q_\lambda$ are weak solutions of the Dirichlet boundary
value problem
\[
\left\{\begin{array}{ll}
l(x,D) u=0\\
j^{m-1}u(0)=0= j^{m-1}u(1).
\end{array}\right.
\]
Let $\mathscr O =  \left \{  z :=\lambda+is \in \C  : 0 < \lambda
<1,  -1 < s <  1 \right \}$. For each $z \in \overline{\mathscr O}$
let us consider the closed unbounded Fredholm operator $A_z$ on
$L^2(J, \C^{mn})$ with domain $\mathscr D( A_z) =\{u \in \mathscr
H^{2m}:j^{m-1} u(0)=0=j^{m-1} u(1) \}$ defined by $ A_zu :=
l_\lambda(x,D)u +isu$. Since $A_{\bar z}= \bar A_z$, it follows that
$A_z$ has a continuous inverse for $s =\Im(z)\neq 0$. Let us
consider the splitting $\C^{2mn}:=(\C^n)^m \times (\C^n)^m$. For
each $w=(w_1, \dots, w_m)\in (\C^n)^m$, let $u_z(w, x)$ be the
unique solution of the Cauchy problem
\[
\left\{\begin{array}{ll}
l(x,D) u+ is u=0\\
u(0)=u'(0)= \dots u^{m-1}(0)=0, \quad u^m(0)=w_1, \dots
,u^{2m-1}(0)= w_m,
\end{array}\right.
\]
and let $\mathscr R_z: \C^{nm}\to \C^{nm}$ be the endomorphism
defined by $\mathscr R_z(w):= j^{m-1}(u_z(w,1))$. Clearly $u \in
\ker A_z$ if and only if $(u^m(0), \dots, u^{2m-1}(0))\in \ker
\mathscr R_z$. From this and by using regularity of the weak
solutions, we get that
the following three statements are equivalent: \\
(i) $\ker  A_{z}\neq \{0\}$; (ii) $ \Im (z) =0 $ and $\lambda = \Re
(z)$ is a conjugate instant; (iii) $ \det \, \mathscr R_z\, = \, 0$.
Because of these three equivalent statements, the function
$\rho(z):=\det \mathscr R_z$ does not vanishes on the boundary
$\partial \mathscr O$.
\begin{defin}
The {\em conjugate index\/} $\icon(q_\Omega)$, of $q_\Omega$ is
minus the winding number of $\rho\vert_{\partial \Omega} \to \C-
\{0\}$ or equivalently minus the Brouwer degree, i.e.
$-\deg(\rho,\mathscr O, 0)$.
\end{defin}
\noindent With this said our main results are

\begin{mainthm}\label{thm:genmorseindex}{\em(Generalized Sturm
Oscillation theorem).\/} With the notation above, we have
\[
\icon(q_\Omega) =\ispec(q_\Omega).
\]
\end{mainthm}

\begin{mainthm}\label{thm:gensturmcompindex}{\em(Generalized Sturm
comparison theorem).} If $\Omega_0, \Omega_1$ are two derivative
dependent Hermitian forms with $ \Omega_0[u](x)\,\leq\,
\Omega_1[u](x)$ for all $x \in J,$ then we have
\[
\ispec(q_{\Omega_0}) \leq \ispec(q_{\Omega_1}).
\]
\end{mainthm}

\begin{section}{Proofs}\label{sec:Proofs}

\proof{{\em(of Theorem \ref{thm:genmorseindex})\/}}. We split the
proof of this result into some steps.\\ \noindent {\em First
step.\/} Notice that the family $A_z$ has the form  $A+ B_z$ where
$A$ is a fixed unbounded closed operator and $B_z$ is an $A$-bounded
perturbation of $A$ depending on $z$ and that the operator valued
one-form $dA_z A_z^{-1}$ is well-defined on $\partial \mathscr O$.
Exactly as in Proposition 4.4 in \cite{MusPejPor03} one proves the
following Lemma.
\begin{lem}\label{thm:prop5.1}
The form   $d{ A}_{z}\, { A}^{-1}_{z} $ takes values in operators of
the trace class and
\begin{equation}\label{eq:tr1}
\icon(q_\Omega) =  -\frac{1}{ 2\pi i}\, \int_{\partial {\mathscr O}}
{\rm Tr}\,d{ A}_{z}\, { A}^{-1}_{z} .
\end{equation}
\end{lem}
\noindent Let us assume now that the path $\{A_\lambda\}_{\lambda
\in J}$ has only regular crossing points with the variety of
singular operators. This means that the quadratic form $\Gamma(A,
\lambda)$ given by the restriction to the kernel of $A_\lambda$ of
$\langle A'_\lambda \cdot,\cdot\rangle$ is nondegenerate. Since
regular crossing points are isolated, there are only finite number
of them. In what follows we will show that if $\mathscr D_j$ is a
small enough neighborhood of a crossing point $\lambda_j \in
\mathscr O$, then
\begin{equation}\label{1}
\frac{1}{2 \pi i}\, \int_{\partial {\mathscr D}_j} {\rm Tr}\, d
A_{z}\, A_{z}^{-1}\, = \, \sgn \,\Gamma ( A , \lambda_j ).
\end{equation}
\noindent For a fixed $j$,  choose a positive number $\mu> 0$ such
that the only point in the spectrum of $ A_{\lambda_j}$ in the
interval $[-\mu , \mu ]$ is 0 and then  choose $\eta$ small enough
such that neither $\mu$ nor $- \mu$ lies in the spectrum of $
A_{\lambda }$ for $ |\lambda-\lambda_j|<\eta.$ For such a $\lambda$,
let $P_{\lambda}$ be the orthogonal projection in $H$ onto the
spectral subspace associated to the part of the spectrum of $
A_{\lambda}$ lying in the interval $[-\mu , \mu ].$ Then $
A_{\lambda}P_\lambda = P_\lambda  A_{\lambda}$ on the domain of
$A_{\lambda}.$ By \cite[Chapter II, Section 6]{Kat80} there exist a
smooth path $ U$ of unitary operators of $H$  defined in $
[\lambda_j -\eta, \lambda_j+\eta] $, such that $ U_{\lambda_j} =
{\rm Id}_H$ and such that
$P_{\lambda } U_\lambda =  U_\lambda  P_{\lambda_j }.$ 
Let us consider the smooth operator valued function $ N_{z }  =
U_\lambda^{-1}\,  A_{z}U_\lambda$ defined on some open neighborhood
of ${\mathscr D}_ j = [\lambda_j -\eta, \lambda_j
+\eta]\times[-1,1]$ together with the differential one-form $\theta
\, = \, d N_{z }\, N_{z }^{-1}$. Clearly $ \theta $  is in the trace
class  and ${\rm Tr} \, d A_z \,  A_z^{-1} \, ={\rm Tr}\, d
N_{z }\, N_{z }^{-1} .$\\
{\em Second step. If $ H_j = \Imm P_{\lambda_j} = \ker
A_{\lambda_j}$, under the splitting  $ H = H_j \oplus H_j ^\perp$,
the one-form $\theta$ splits into $ \theta_0 = d N_{z } N_z^{-1}
|_{H_j} \ \textrm{and} \ \theta_1 = d N_z N_z^{-1}\,|_{H_j^\perp} $.
Thus by taking traces we have
\begin{equation}\label{eq:6.8jac1}
\frac{1}{ 2\pi i } \,\int_{\partial {\mathscr D}_j} \,    {\rm Tr}
\, \theta \, = \, \frac{1}{ 2\pi  i} \, \int_{\partial {\mathscr
D}_j} \,{\rm Tr}\, \theta_0 \,+\, \frac{1}{2\pi i  } \,
\int_{\partial {\mathscr D}_j} \, {\rm Tr} \, \theta_1
\end{equation}
where the last term in \eqref{eq:6.8jac1} vanishes.\/} \\ \noindent
In fact, $N_z\vert_{H_j^\perp}$ is invertible on $\overline{\mathscr
O}$ and hence the one form $\theta_1$ is exact. Now, let us consider
the path of symmetric endomorphisms $M \colon [\lambda_j -\eta,
\lambda_j + \eta ] \to {\mathscr L} (H_j) $ given by $ M_\lambda
={N_\lambda}\vert_{H_j}$. By elementary calculations
it follows that 
$\sgn \,\Gamma(M_\lambda, \lambda_j) = \sgn \, \Gamma( A,
\lambda_j).$\\
{\em Third step. If $M$ is the path defined above  and if $ M_z =
M_\lambda +is \Id,$ for $ z \in {\mathscr D}_j$, then
\[
\sgn \,\Gamma(M, \lambda_j)   = \,\frac{1}{2\pi i} \int_{\partial
{\mathcal D}_j} \,{\rm Tr}\, dM_z M_z^{-1} .
\] \/}\noindent
Using  \cite[Chapter II, Theorem 6.8]{Kat80}, one reduces
$M_\lambda$ to a diagonal form and the result follows by direct
integration. This together with \eqref{eq:6.8jac1} proves \eqref{1}. \\
\noindent {\em Fourth step.\/} Using regularity it is easy to see
that any regular crossing point for the path $A$ is also a regular
crossing point for the path $q$ i.e., the crossing form defined as
the restriction $\Gamma(q, \lambda)$ of $\dot q_\lambda$ to $\ker
q_\lambda$ is non-degenerate. Moreover the crossing forms $\Gamma(A,
\lambda)$ and $\Gamma(q, \lambda)$ coincide. Since for path with
only regular crossings the spectral flow is given by the sum of the
signatures of the crossing forms it follows from Lemma
\ref{thm:prop5.1} and formula \eqref{1} that theorem
\ref{thm:genmorseindex} holds for paths with only regular crossings.
In order to conclude remains to show that it is possible to extend
the above calculation to general paths having only regular
crossings. To do so we will apply a perturbation argument of Robbin
and Salamon to the path of operators obtained by the restriction of
the family $A_z$ to the real line. By Theorem 4.22
in\cite{RobSal95}, we can find a $\delta>0$ such that
$A_\lambda^\delta:= A_\lambda + \delta \Id$ is a path of
self-adjoint Fredholm operators with only regular crossing points.
Let $ \mathscr R^\delta_{z}$ be the matrix associated to the
perturbed family $ A^{\delta}_{z }$ and let $ \rho^\delta (z) :=
\det\, \mathscr R^\delta_{z}$. Taking closed disjoint neighborhoods
${\mathscr D}_j$ of $(\lambda_j,0)$ in $ {\mathscr O}$ and summing
over all crossing points, we obtain
\begin{eqnarray*}
\icon(q_\Omega) \, = \,-\sum_{j=1}^{k}\deg (  \rho^\delta ,
{\mathscr D}_j , 0 ) = \sum_{j=1}^{k} \,-\frac{1}{2 \pi i}\,
\int_{\partial {\mathscr D}_j}  {\rm Tr}\,d A^\delta_{z} (
A^\delta_z)^{-1}= -\sum_{j=1}^{k} \sgn \,\Gamma ( A , \lambda_j )=
-\spfl (q_\lambda, J)
\end{eqnarray*}
by the homotopy invariance of the spectral flow. Now Theorem
\ref{thm:genmorseindex} is proved.  \finedim

\proof{{\em of Theorem \ref{thm:gensturmcompindex}}}. Since the top
order terms coincide, the difference between the two Fredholm
Hermitian forms $q^1$ and $q^2$ is weakly semi-continuous. Thus in
order to conclude it is enough to show that if $q^1$ and $q^2$ are
two admissible families of Fredholm Hermitian forms on $\mathscr
H_0^m$ whose difference is weakly semi-continuous and such that
$q^1_0 = q^2_0,$ $q^2_1\leq q^1_1,$ then $\spfl(q^2) \leq
\spfl(q^1).$ The rest of the proof is devoted to show this.

It is easy to see that there exists a small $\eta >0$ such that
$\tilde q^1_\lambda = q^1_\lambda + \lambda  \eta ||\cdot ||^2 $ is
a  family of Fredholm Hermitian forms that is homotopic to $q^1$ by
$ (t,\lambda) \to q^1_\la + t \lambda \eta ||\cdot ||^2$  and
therefore $\spfl(\tilde q^1 ) = \spfl(q^1)$. But now $ \tilde q^1_0
= q^2_0 $ and $ q^2_1 < \tilde q^1_1$. Given the family of  Fredholm
Hermitian forms defined on $ T$ by $\phi (s,\lambda ) = s \tilde
q^1_\lambda + (1-s) q^2_\lambda$ and by using the free homotopy
property of the spectral flow for closed paths, it follows that the
restriction of $\phi$ to $  \partial T$ has spectral flow zero.
Furthermore, $\spfl (\widetilde q^1  ) - \spfl (q^2) = \spfl (\rho),
$ where $\rho(s)= (1-s) q^2_1 + s \widetilde q^1_1.$ Since $\dot
\rho = \widetilde q^1_1-q^2_1$ is positive definite at each crossing
point, all crossing points are regular, and  each gives a positive
contribute $\dim\ker \ \rho(s)$ to the spectral flow of $\rho.$ Thus
$ \spfl (\rho) \geq 0$ and hence $\spfl (q^2) \leq \spfl (q^1).$
\finedim
\end{section}
\subsection*{Acknowledgements} I would like to thank Jacobo
Pejsachowicz for many useful conversations. I am also very grateful
to the anonymous referee for his numerous suggestions, criticisms
and for having pointed out some crucial mistakes in the original
version.


\end{document}